\documentclass[12pt]{amsart}

\usepackage{latexsym, amssymb, amscd, amsfonts}

\usepackage{mathpple}
\usepackage{pstricks}
\usepackage{pst-poly}
\usepackage{multido}
\usepackage{pst-node}
\usepackage[all]{xypic}
\usepackage{ifthen}

\newboolean{showpicture}
\setboolean{showpicture}{false}

\newboolean{bourbaki}
\setboolean{bourbaki}{true}

\newboolean{draft}
\setboolean{draft}{true}

\newcommand{\cplus}{\hbox{$\subset${\raise0.3ex\hbox{\kern -0.55em${\scriptscriptstyle +}$}}\ }}
\newcommand{\clplus}{\hbox{$\subset${\raise0.3ex\hbox{\kern -0.55em ${\scriptscriptstyle +}$}}\ }}
\newcommand{\crplus}{\hbox{$\supset${\raise1.05pt\hbox{\kern -0.55em ${\scriptscriptstyle +}$}}\ }}

\newcommand{\N}{\mathbf{N}}

\renewcommand{\phi}{\varphi}
\newcommand{\OO}{\mathcal{O}}
\newcommand{\LL}{\mathcal{L}}

\newcommand{\C}{\mathbf{C}}
\newcommand{\Z}{\mathbf{Z}}

\newcommand{\hh}{\mathfrak{h}}
\newcommand{\bb}{\mathfrak{b}}

\newcommand{\refle}[1]{Lemma \ref{#1}}
\newcommand{\refth}[1]{Theorem \ref{#1}}

\newcommand{\refsec}[1]{section~ \ref{#1}}
\newcommand{\refeq}[1]{(\ref{#1})}
\newcommand{\eqdef}{:=}

\def\gb{\mathfrak{b}}

\def\gg{\mathfrak{g}}
\def\gh{\mathfrak{h}}

\def\C{\mathbb{C}}

\def\K{\mathbb{K}}

\def\N{\mathbb{N}}

\def\P{\mathbb{P}}

\def\Z{\mathbb{Z}}

\def\cL{\mathcal{L}}

\def\cO{\mathcal{O}}
\def\cP{\mathcal{P}}

\def\cS{\mathcal{S}}
\def\cT{\mathcal{T}}

\def\inj{\varinjlim}
\def\proj{\varprojlim}

\def\nek{\text{\hbox{$\simeq$ \kern-.95em \hbox{$/$ \kern.05em}}}}
\def\vep{\varepsilon}

\def\cplus{\hbox{$\subset${\raise1.05pt\hbox{\kern -0.55em
${\scriptscriptstyle +}$}}\ }}
\def\bcplus{\hbox{$\supset${\raise1.05pt\hbox{\kern -0.55em
${\scriptscriptstyle +}$}}\ }}

\def\ctimes{\hbox{$\times${\raise1.1pt\hbox{\kern -0.27em
${\scriptscriptstyle |}$}}\ }}

\def\udarrow{\hbox{$\nearrow${\kern -0.97em$\searrow$}\ }}

\def\bctimes{\hbox{$\times${\raise1.1pt\hbox{\kern -.74em
${\scriptscriptstyle |}$}}\ }\,\,}

\newtheorem{theorem}[equation]{Theorem}

\newtheorem{lemma}[equation]{Lemma}

\newtheorem{corollary}[equation]{Corollary}

\newtheorem{proposition}[equation]{Proposition}

\theoremstyle{remark}
\newtheorem{example}[equation]{Example}

\newtheorem{remark}[equation]{Remark}

\title {A Bott--Borel--Weil theorem for diagonal ind--groups}
\author{Ivan Dimitrov}
\address{Department of Mathematics and Statistics, Queen's University, Kingston,
Ontario,  K7L 3N6, Canada} 
\email{dimitrov@mast.queensu.ca} 
\thanks{Research partially supported by an NSERC grant}
\author{Ivan Penkov}
\address{Jacobs University Bremen, Campus Ring 1, 28759 Bremen, Germany} 
\email{i.penkov@jacobs-university.de} 
\thanks{Research partially supported by DFG Grant PE980/2-1}

\begin{document}
\maketitle
\begin{abstract}
A diagonal ind--group is a direct limit of classical affine algebraic groups of growing rank under a class of
embeddings which contains the embedding
$$
SL(n)\to SL(2n), \quad \quad
M\mapsto \left(\begin{array}{cc}M & 0\\ 0 & M \end{array}\right)
$$
as a typical special case. If $G$ is a diagonal ind--group and $B\subset G$ is a Borel ind--subgroup, 
we consider the ind--variety $G/B$ and compute the cohomology $H^\ell(G/B,\mathcal{O}_{-\lambda})$ 
of any $G$--equivariant line bundle $\mathcal{O}_{-\lambda}$ on $G/B$. It has been known that, for a generic $\lambda$, 
all cohomology groups of $\cO_{-\lambda}$ vanish, and that a non--generic equivariant 
line bundle $\cO_{-\lambda}$ has at most one 
nonzero cohomology group. The new result of the present paper is a precise description of when 
$H^j(G/B,\mathcal{O}_{-\lambda})$ is nonzero and the proof of the fact that, whenever nonzero, 
$H^j(G/B, \mathcal{O}_{-\lambda})$ is a $G$--module dual to a highest weight module. 
The main difficulty is in defining an appropriate analog $W_B$ of the Weyl group, so that the action of $W_B$
on weights of $G$ is compatible with the analog of the Demazure "action" of the Weyl group on the cohomology 
of line bundles. The highest weight corresponding to $H^j(G/B, \cO_{-\lambda})$ is then computed
by a procedure similar to that in the classical Bott--Borel--Weil theorem.
\end{abstract}

\section*{Introduction}
The classical Bott--Borel--Weil theorem is a cornerstone of geometric representation theory. 
In the late 1990's Joseph A. Wolf and his collaborators became interested in extending the theorem 
to direct limit Lie groups, and since then have made essential progress, see \cite{NRW}, \cite{W}. 
In the context of direct limit algebraic groups, i.e. ind--groups, the problem has been addressed in 
our joint paper \cite{DPW}. In that paper a quite general theorem has been proved (concerning 
infinite--rank equivariant bundles on locally proper homogenous ind--varieties), under the condition 
that the ind--group considered is root reductive, the definition see in \refsec{sec21} below. The known results become much 
sketchier  when this condition is dropped. The purpose of the present paper is to consider in detail 
the most interesting class of ind--groups beyond the root reductive ones, that of diagonal ind--groups.

Recall that a locally affine ind--group G is the direct limit of embeddings of connected affine algebraic groups
$$
G_{1}\rightarrow G_{2}\rightarrow\ldots \quad .
$$
The Bott--Borel--Weil paradigm for ind--groups is concerned with the computation of the cohomology 
of a $G$--equivariant line bundle $\mathcal{O}_{-\mu}$ on the ind--variety $G/B$, where 
$B=\lim\limits_{\rightarrow}B_{n}$ is the direct limit of Borel subgroups $B_{n}\subset G_{n}$ 
with $B_{n-1} = B_{n} \cap G_{n-1}$. In the classical case $G$ is a connected affine algebraic group
and the result (due to Borel--Weil \cite{BW} and Bott \cite{B}, see also \cite{D1}, 
\cite{D2}) is that the simple (finite--dimensional) $G$--module $V_{B}(\lambda)^*$ with $B$--highest 
weight $\lambda$ occurs as the unique nonzero cohomology group of each of the sheaves 
$\mathcal{O}_{-w\cdot \lambda}$, where $w$ runs over the Weyl group W and $\cdot$ stands for 
the "dot action" of $w$ on $\lambda$. More precisely, $V_{B}(\lambda)^*$ occurs as the
cohomology group of $\mathcal{O}_{-w\cdot \lambda}$ in degree $\ell(w)$, where $\ell(w)$ is the length of $w$ with 
respect to the simple roots of $B$.

In contrast with this result, in the infinite--dimensional case it is not difficult to see that a generic 
line bundle $\mathcal{O}_{-\mu}$ is acyclic, i.e. all its cohomology groups vanish.  Wolf has 
introduced the condition of cohomological finiteness of a weight $\mu$ (see \cite{W} and 
compare with \cite{DPW}), which is equivalent to the condition that $\mathcal{O}_{-\mu}$ has a 
unique non--vanishing cohomology group. If $\mu$ is dominant, then this cohomology group is 
$H^0(G/B, \mathcal{O}_{-\mu})$, and in this case it is easy to show that 
$H^0(G/B, \mathcal{O}_{-\mu})$ is the (algebraic) dual of the simple $B$--highest weight 
G--module $V_{B}(\mu)$.

What is not known in general is whether all higher cohomology groups 
$H^j(G/B, \mathcal{O}_{-\mu})$ are also dual to $B$--highest weight modules. This problem has 
been open since the late 1990's, and the main result of the present paper is that for any locally 
simple diagonal ind--group $G$ (see the definition in \refsec{sec21}), all nonzero cohomology 
groups $H^j(G/B, \mathcal{O}_{-\lambda})$ are indeed dual to simple $B$--highest weight 
modules. The proof is a mixture of combinatorics and geometry. The most important new 
idea is to consider the intermediate algebraic groups $\widetilde{G}_{n} \cong G_n \times G_n \times \ldots \times G_n$,
$$
G_{n}\to \widetilde{G}_{n} \to G_{n+1},
$$
introduced in \refsec{sec21}. They arise naturally from the diagonal embeddings 
$G_{n} \to G_{n+1}$. The corresponding homogenous spaces 
$\widetilde{G}_{n} / \widetilde{B}_{n}$, where 
$\widetilde{B}_{n} = B_{n+1} \cap \widetilde{G}_{n}$, play a key role in the proof. 
More precisely, the realization of $\mathcal{O}_{-\lambda}$ as a line bundle on both $G/B$ 
and on $\displaystyle\lim_{\to} \tilde{G}_n/\tilde{B}_n$ enables us to reduce the problem 
of studying the cohomologies $H^j(G/B, \mathcal{O}_{-\lambda})$ to two finite--dimensional
problems --- one concerns the embeddings $G_{n}\to \widetilde{G}_{n}$ and the one concerns
the embeddings $\widetilde{G}_{n} \to G_{n+1}$. For the second problem we use a recent result of  
Valdemar Tsanov, 
which allows us to obtain a strong condition on the weight $\lambda$ so that 
$H^j(G/B, \mathcal{O}_{-\lambda}) \neq 0$; under this condition we then apply a result of Mike Roth and 
the first named author to the embedding 
$G_{n}\to \widetilde{G}_{n}$.
The final result, \refth{th2}, 
is absolutely similar to the classical Bott--Borel--Weil Theorem with the only exception that the 
"Weyl group" $W_B$, relevant for $G/B$,  depends on the choice of Borel ind--subgroup $B$. 

\vskip.5cm

\noindent
\textbf{Acknowledgement. }We thank the Mathematisches Forschungsinstitut Oberwolfach, 
Germany, and the Banff International Research Station, Canada, where parts 
of this work were done. I. D. acknowledges the hospitality of the ICTS at Jacobs University Bremen.

\vskip.5cm

\section{Diagonal ind--groups: definitions and notation}\label{sec21}
We work over an algebraically closed field $\K$ of characteristic $0$.
If $V$ is a vector space, we set $V^{\oplus k}=\underbrace{V\oplus\dots\oplus V}_{k \mathrm{~times}}$.

Throughout this paper, {\it classical group} will be an abbreviation for a connected (affine)
algebraic group $G$ whose Lie algebra is a simple classical Lie algebra. An embedding 
$G \to G'$ of classical groups is ${diagonal}$ if the induced injection of Lie 
algebras $\frak{g} \to \frak{g'}$ has the following property: the natural 
representation of $\frak{g'}$ considered as a $\frak{g}$--module is isomorphic to a direct 
sum of copies of the natural representation of $\frak{g}$, of its dual, and of the trivial 
representation. (If $\frak{g}=so, sp$ the natural representation is self--dual, hence in this 
case the natural representation of $\frak{g}'$ must simply be a direct sum of copies of the 
natural and trivial representations). 
If $\gg$ and $\gg'$ are reductive Lie algebras, an injective Lie algebra homomorphism
 $\gg\to\gg'$ is a {\it root injection}, if for any Cartan subalgebra $\hh\subset\gg$ there exists 
 a Cartan subalgebra $\gh' \subset \gg'$
 containing the image of $\gh$ and such that any $\hh$--root space of $\gg$ is mapped to 
 precisely one $\hh'$--root space of $\gg'$.
 An embedding $G \to G'$ of reductive affine algebraic groups is a {\it root embedding} if the 
 corresponding injection $\gg \to \gg'$ is a root injection.

By definition a {\it diagonal ind--group} $G$ as the direct limit of a sequence of diagonal 
embeddings of classical groups
\begin{equation}\label{eq11}
G_{1} \rightarrow \ldots \rightarrow G_{n}\rightarrow G_{n+1}  \rightarrow \ldots \quad .
\end{equation}
The group $G$ is called {\it pure} if, for large enough $n$, the natural representation of $\gg_{n+1}$ contains no
trivial $\gg_n$--constituents.

\vskip.5cm

\begin{example} \label{ex1}
A diagonal embedding of classical groups $G \to G'$ of type $A$ can be realized 
in  matrix form as
$$ M \longmapsto
\begin{pmatrix} M  & & &  & & & & \\
  &  \ddots      &  &  & & & &  \\
  &              & M &  & & & & \\
  &              &   & (M^\intercal)^{-1} & & & &\\
  & & & & \ddots & & & \\
  & & & & & (M^{\intercal})^{-1} & & \\
  & & & & & & 1 & \\
  & & & & & & & \ddots & \\
  & & & & & & & & 1 \\
\end{pmatrix}
$$ 
with $k$ copies of $M$, $l$ copies of $(M^{\intercal})^{-1}$, and $t$ copies of the one--by--one matrix with entry one. 
Therefore, any diagonal ind--group 
of type $A$ is obtained by
iterating such embeddings with varying parameters $k, l$, and $t$. In particular, the 
ind--group $SL(\infty)$ can be defined  as a  diagonal ind--group of type $A$ with
$k= t =1$, $l=0$ at each step. To define the diagonal ind--group $SL(2^\infty)$ we set  $G_1 := SL(2)$ and then put  
$k=2$, $l=t=0$ at each step. It is easy to check that, up to isomorphism, $SL(\infty)$ does not depend on the choice of $n_1$ where
$G_1 = SL(n_1)$. 
The ind--group $SL(2^\infty)$ is pure while $SL(\infty)$ is not.
\hfill $\square$
\end{example}

\vskip.5cm

In this paper we consider two types of $G$--modules, defined respectively as direct or inverse limits of 
finite--dimensional $G_{n}$--modules.
Fix $G = \lim\limits_{\rightarrow}G_{n}$, and let
\begin{equation}\label{eqGmod}
    V_{1} \rightarrow \ldots  \rightarrow V_{n} \rightarrow V_{n+1} \rightarrow \ldots ,
\end{equation}
(respectively,
\begin{equation}\label{eq313}
\ldots \rightarrow Y_{n+1} \rightarrow Y_{n} \rightarrow \ldots \rightarrow Y_1  \mathrm{)}
\end{equation}
be a direct (resp. inverse) system of finite--dimensional $G_{n}$--modules. By a {\it $G$--module} will mean 
the direct limit of a system \refeq{eqGmod} endowed with $G_{n}$--module structures for all $n$, up to an 
isomorphism, and by a {\it dual $G$--module}  mean the projective limit of a system \refeq{eq313} 
endowed with $G_{n}$--module structures. It is clear that if 
$V =\lim\limits_{\rightarrow}V_{n}$ is a $G$--module, then $V^{*} =\lim\limits_{\leftarrow}V_{n}^{*}$ is a 
dual $G$--module. Conversely, if $Y = \lim\limits_{\leftarrow}Y_{n}$ is a dual $G$--module, then 
$\lim\limits_{\rightarrow}Y_{n}^{*}$ is a $G$--module.

For the rest of the paper we fix an exhaustion $G = \inj G_n$ of $G$ by simply--connected classical groups
of the same type $A, B, C$, or $D$. In particular, every direct system (\ref{eq11}) we consider has 
a well--defined type.
In general $G$ may have exhaustions of different type, however we will use the term "type of $G$" to 
refer to the type of the fixed exhaustion. The corresponding exhaustion of $\gg$ is then $\gg = \inj \gg_n$.
We denote  the rank of $\gg_n$ by $r_n$.

For the purposes of this paper, we define a {\it Cartan subgroup} $H$ of $G$ as a direct limit 
of Cartan subgroups $H_{n} \subset G_{n}$. 
The corresponding Lie algebra $\frak{h}$ is then the direct limit of Cartan subalgebras 
$\frak{h}_{n} \subset \frak{g}_{n}$ such that $\frak{h}_{n} = \frak{h}_{n+1} \cap \frak{g}_{n}$. We 
fix once and for all a Cartan subgroup $H = \inj H_n$ of $G$ with corresponding Cartan subalgebra
$\gh = \inj \gh_n$ of $\gg$. The weights of $\gg_n$ are expressed in terms of the standard functions 
$\vep_n^1, \ldots, \vep_n^{r_n + 1} \subset \gh_n^*$ if $G$ is of type $A$ or 
$\vep_n^1, \ldots, \vep_n^{r_n} \subset \gh_n^*$ otherwise. These functions are determined by the choice of the Cartan
subalgebra $\gh_n \subset \gg_n$. The weights of the natural representation of $\gg_n$ are as follows: for $G$ of type $A$ 
they are $\vep_n^1, \ldots, \vep_n^{r_n + 1}$; for $G$ of type $B$ --- $\pm \vep_n^1, \ldots, \pm \vep_n^{r_n}, 0$; and
for $G$ of type $C$ or $D$ --- $\pm \vep_n^1, \ldots, \pm \vep_n^{r_n}$.  Since $\gh_n \subset \gh_{n+1}$, the 
$\gh_{n+1}$--weight spaces of the natural representation of $\gg_{n+1}$ restrict to $\gh_n$--weight spaces. In particular,
 $\vep_{n+1}^i$ restricts to $\pm\vep_n^j$ for some $j$, or to $0$. 

Denote the injection $\gg_n \to \gg_{n+1}$ by $\delta_n$. 
We will now define a subalgebra 
$\widetilde{\gg}_{n} \cong \gg_n^{\oplus s_{n}}$
of $\gg_{n+1}$, where $s_{n}$ is the 
total multiplicity of all  nontrivial simple constituents of the natural representation of $\gg_{n+1}$
considered as a $\gg_n$--module. 
Note first that the $\gh_{n+1}$--weight decomposition of the natural representation of $\gg_{n+1}$ determines a 
unique decomposition of each nontrivial isotypic $\gg_n$--component 
as a direct sum of simple constituents. To define the subalgebra $\tilde{\gg}_{n}$ it suffices to define its simple ideals:
if $G$ is of type $A$, each simple ideal of $\tilde{\gg}_{n}$ equals the traceless endomorphisms of a simple 
nontrivial constituent of the natural representation of $\gg_{n+1}$; if $G$ is of type $B, C$, or $D$, each simple ideal of 
$\tilde{\gg}_{n}$ is the Lie algebra of orthogonal or respectively symplectic endomorphisms of a simple nontrivial constituent
 of the natural representation of $\gg_{n+1}$. In all cases, there is an obvious injective homomorphism 
$\varphi_{n}: \gg_n \to \tilde{\gg}_{n}$ such that the diagram 
\begin{equation} \label{eq12}
\xymatrix{
\gg_n  \ar[rr]^{\delta_{n}} \ar[dr]_{\varphi_{n}}  & &\gg_{n+1} \\
& \tilde{\gg}_{n}
\ar[ur]_{\kappa_{n}} &}
\end{equation}
is commutative, $\kappa_{n}$ being the inclusion. Moreover,
if $G$ is of type $A, C$, or $D$, the map $\kappa_{n}$ is a root injection. 

If $G$ is of type $B$ and $s_n >1$, $\kappa_{n}$ is no longer a root injection, however
we can still factor $\kappa_{n}$ as
$$
\widetilde{\gg}_{n} \stackrel{\psi_{n}}{\longrightarrow}  \bar{\gg}_{n} \stackrel{\theta_{n}}{\longrightarrow}  \gg_{n+1} 
$$
so that $\theta_{n}$ is a root injection and $\psi_{n}$ is "close" to a root injection. To construct this factorization,
recall that  $\gg_n \cong B_{r_n}$, $\widetilde{\gg}_{n} \cong B_{r_n}^{\oplus {s_n}}$, 
$\gg_{n+1} \cong B_{r_{n+1}}$, and  that the natural representation $V_{n+1}$ of $\gg_{n+1}$
when considered as a $\gg_n$--module contains $s_n$ copies of the natural representation of $\gg_n$ and $z_n$ copies of the trivial
representation. Note that $2r_{n+1}+1 = s_n(2r_n+1) + z_n$, hence $s_n$ and $z_n$ are distinct modulo $2$. The smallest
interesting case is when $s_n=2$ and $z_n=1$. In this case the $\gg_n$--module decomposition of the
natural representation of $\gg_{n+1}$ is $V_{n+1} = V^1 \oplus V^2 \oplus \K$, where $V^1$ and $V^2$ are the two copies of the natural
representation of $\gg_n$. We set $\bar{\gg}_n := so(V^1 \oplus V^2) \cong D_{2r_n+1}$ and consider the natural
injections $\widetilde{\gg}_n  \stackrel{\psi_n}{\to} \bar{\gg}_n  \stackrel{\theta_n}{\to} \gg_{n+1}$. 
The assumption that $\gh_n$ is contained in $\gh_{n+1}$ 
 ensures that $\gh_{n+1}$ is contained in $\bar{\gg}_n$
and, consecutively, $\theta_n$ is a root injection. Furthermore, for each long root of $\widetilde{\gg}_n$,
the corresponding root space is mapped via $\psi_n$ into a root space of $\bar{\gg}_n$. 
In the case when $s_n$ and $z_n$ are arbitrary,
i.e. when $V' = V^1 \oplus \ldots \oplus V^{s_n} \oplus \K^{z_n}$, we combine $V^1, \ldots, V^{s_n}$ into pairs when 
$s_n$ is even, and into pairs and a
single element when $s_n$ is odd, and set 
$\bar{\gg}_n := so(V^1 \oplus V^2) \oplus \ldots \oplus so(V^{s_n-1} \oplus V^{s_n})$ in the former case
and $\bar{\gg}_n := so(V^1 \oplus V^2) \oplus \ldots \oplus so(V^{s_n-2} \oplus V^{s_n-1}) \oplus so(V^{s_n})$ in the latter. 
The injections $\widetilde{\gg}_n  \stackrel{\psi_n}{\to} \bar{\gg}_n  \stackrel{\theta_n}{\to} \gg_{n+1}$ are
defined in the obvious way, $\theta_n$ is a root injection, and $\psi_n$ maps roots spaces corresponding to long roots of
$\widetilde{\gg}_n$ to root spaces of $\bar{\gg}_n$. As a result of this construction, we obtain a refinement of
diagram (\ref{eq12}) as follows:
\begin{equation} \label{eq13}
\xymatrix{
\gg_n  \ar[rrrr]^{\delta_{n}} \ar[dr]_{\varphi_{n}}  & && &\gg_{n+1} \\
& \widetilde{\gg}_{n} \ar[rr]_{\psi_{n}} \ar[rrru]^{\kappa_{n}}&& \bar{\gg}_{n}
\ar[ur]_{\theta_{n}} &}.
\end{equation}

\vskip.5cm

\begin{remark} \label{remark5}
Note that the subalgebra $\bar{\gg}_n$ above depends on the way we combine $V^1, \ldots, V^{s_n}$ into pairs. 
In the proof of Corollary \ref{cor51} below we consider diagrams analogous to \eqref{eq13}  for different choices for 
$\bar{\gg}_n$ in the case when $s_n$ is odd and greater than one.
\end{remark}

\vskip.5cm

 \section{Borel subalgebras, dominant weights, and highest weight modules}
For the purposes of this paper we adopt the following definition of a Borel subgroup: {\it $B$ is a Borel subgroup of $G$} if
$B = \inj B_n$, where $B_n$ is a Borel subgroup of $G_n$ for every $n$. All Borel sugroups we consider contain
the fixed Cartan subgroup $H$. 
The corresponding Lie algebra $\gb$  then contains $\gh$ and is the 
direct limit of Borel subalgebras $\gb_{n} \subset \frak{g}_{n}$ containing the fixed Cartan subalgebras $\gh_n$.
Note that we have $\gb_n = \gb_{n+1} \cap \gg_n$.  

The Borel subalgebras of $\gg_n$ which contain $\gh_n$ correspond to linear orders on the weights of the natural representation of $\gg_n$.
More precisely (see also \cite{DP}), the Borel subalgebras of $\gg_n$ correspond to
\begin{itemize}
\item the linear orders on the set $\{\vep_n^1, \ldots, \vep_n^{r_n+1}\}$ for type $A$;
\item the linear orders compatible with multiplication by $-1$ on the set $\{\pm \vep_n^1, \ldots, \pm \vep_n^{r_n}, 0\}$ for 
type $B$;
\item the linear orders compatible with multiplication by $-1$ on the set $\{\pm \vep_n^1, \ldots, \pm \vep_n^{r_n}\}$ for
type $C$;
\item the linear orders compatible with multiplication by $-1$ on the set $\{\pm \vep_n^1, \ldots, \pm \vep_n^{r_n}\}$ for type 
$D$.
\end{itemize}
Here "compatible with multiplication by $-1$" means that $\vep_n^i < \pm \vep_n^j$ is equivalent to $\mp \vep_n^j < - \vep_n^i$.
The above correspondence is a bijection in types $A, B$, and $C$; in type $D$ each Borel subalgebra corresponds to exactly
two orders as above since the smallest element 
$\pm \vep_n^i$ such that $\pm \vep_n^i > \mp \vep_n^i$ can be interchanged with its
opposite without changing the Borel subalgebra.

The condition $\gb_n = \gb_{n+1} \cap \gg_n$ is equivalent to the fact that the order on the weights of the natural 
representation of $\gg_{n+1}$ restricts to the order (or one of the two orders in type $D$)
on the weights of the natural representation of $\gg_{n}$. In this way we can say that a Borel subalgebra $\gb = \inj \gb_n$
is determined by a projective system of linear orders on the weights of the natural representations of $\gg_n$. 
Note that in type $A$ the weights of $\gg_{n+1}$
corresponding to constituents isomorphic to the dual of the natural representation of $\gg_n$ restrict to $-\vep_n^i$.

\vskip.5cm

\begin{example} \label{ex2}
(i) Let  $G = SL(2^\infty)$. Then $\gg_n = sl(2^n)$ with weights $\{\vep_n^1, \ldots, \vep_n^{2^n}\}$ 
of the natural representation.  
The weights $\vep_{n+1}^i$ and $\vep_{n+1}^{2^n + i}$ restrict to $\vep_n^i$ for
$1 \leq i \leq 2^n$. The projective system of orders
$$
\vep_n^1 > \vep_n^2 > \ldots > \vep_n^{2^n}
$$
defines the Borel subgroup of $G$ consisting of upper triangular matrices in the realization of $G$ from Example \ref{ex1}.
We will call this Borel subgroup {\it the upper triangular Borel subgroup} of $SL(2^\infty)$.  

\noindent
(ii) A more interesting example of a Borel
sugroup of $SL(2^n)$ is provided by the projective systems of orders
$$
\vep_n^1 > \vep_n^{2^{n-1}+1} > \vep_n^2 > \vep_n^{2^{n-1}+2} > \ldots > \vep_n^{2^{n-1}} > \vep_n^{2^n}.
$$
We will call this Borel subgroup {\it the interlacing Borel subgroup} of $SL(2^\infty)$.

\noindent
(iii) Let $G_n = Sp(2(2^n-1))$ and let the embedding $G_n \to G_{n+1}$ be determined by the condition that
the natural representation of $\gg_{n+1}$ contains two copies of the natural representation of $\gg_n$ and two
copies of the trivial representation. The resulting ind--group is not pure; we denote it by $Sp(2^\infty +1)$. 
The weights of the natural representation of $\gg_n$ are $\pm \vep_n^1, \ldots,
\pm \vep_n^{2^n-1}$. We assume that $\vep_{n+1}^{1+i}$ and $\vep_{n+1}^{1 + 2^n + i}$ restrict to $\vep_n^{1+i}$
for $1 \leq i \leq 2^n-2$, while $\vep_{n+1}^1$ restricts to $0$. The projective system of orders
$$
\begin{array}{rl}
 \vep_n^1 > & \vep_n^2 > \vep_n^{2 + 2^{n-1}} > \ldots > \vep_n^{2^{n-1}-1} > \vep_n^{2^n-1} \\
> & - \vep_n^{2^n-1} > - \vep_n^{2^{n-1}-1} > \ldots > -  \vep_n^{2 + 2^{n-1}} > - \vep_n^2 > - \vep_n^1 
\end{array} 
$$
defines an {\it interlacing Borel subgroup $B$} of $Sp(2^\infty+1)$.
\hfill $\square$
\end{example}

\vskip.5cm

A {\it weight $\lambda$ of $G$} is by definition an inverse system of weights of $G_n$, i.e. 
a sequence $\{\lambda_n\}$ of integral weights of $\gg_n$ such that $\lambda_{n+1}$
restricts to $\lambda_n$ for every $n$. We use the notation $\lambda = \proj \lambda_n$ to indicate that the sequence
$\{\lambda_n\}$ defines the weight $\lambda$.  $\cP$ stands for the set of weights of $G$. As in the finite--dimensional case,
for every Borel subgroup $B \subset G$, $\cP$ is in a natural bijection with the one--dimensional $B$--modules.
A weight $\lambda \in \cP$ is {\it $B$--dominant} (or, simply, {\it dominant} if $B$ is clear from the context) if $\lambda_n$
is a $B_n$--dominant weight for every $n$; the set of $B$--dominant weights will be denoted by $\cP_B^+$ (respectively, by
$\cP^+$). The fundamental $\gb_n$--weights of $\gg_n$ (in the standard order on the nodes of the Dynkin diagram of
$\gg_n$) will be denoted by $\omega_n^1, \ldots, \omega_n^{r_n}$. 

\vskip.5cm

\begin{example} \label{ex3}
We discuss $\cP_B^+$ for each of the Borel subgroups from Example \ref{ex2}.

\noindent
(i) Consider $\lambda = \proj \lambda_n \in \cP$. Let 
\begin{equation} \label{eq31}
\lambda_n  =  \lambda_n^1 \vep_n^1 + \ldots + \lambda_n^{2^n} \vep_n^{2^n} 
=  a_n^1 \omega_n^1 + \ldots + a_n^{2^n-1} \omega_n^{2^n-1}.
\end{equation}
Since $a_n^i = \lambda_n^i - \lambda_n^{i+1}$ for $1 \leq i \leq 2^n-1$, 
the fact that $\lambda_{n+1}$ restricts to $\lambda_n$ is 
equivalent to the equations
\begin{equation} \label{eq32}
a_n^1 = a_{n+1}^1 + a_{n+1}^{2^n+1}, a_n^2 = a_{n+1}^2 + a_{n+1}^{2^n+2}, \ldots, 
a_n^{2^n-1} = a_{n+1}^{2^n-1} + a_{n+1}^{2^{n+1}-1},
\end{equation}
and $\lambda \in \cP^+$ is equivalent to the condition that $a_n^i \in \Z_{\geq 0}$ for every $n$ and every $1 \leq i \leq 2^n-1$.
As \eqref{eq32} shows, every $B_n$--dominant weight $\lambda_n$ of $G_n$ is  the restriction of infinitely many 
$B_{n+1}$--dominant weights of $G_{n+1}$. More precisely, there are finitely many choices for the parameters
$a_{n+1}^1, \ldots, a_{n+1}^{2^n-1}, a_{n+1}^{2^n+1}, \ldots, a_{n+1}^{2^{n+1}-1}$,  and the parameter 
$a_{n+1}^{2^n}$ can be chosen as any element of $\Z_{\geq0}$. 
In particular, $\cP^+$ is not finitely generated and contains the lattice points of 
an open $n$--dimensional cone for every $n$.

\noindent
(ii) As in (i) above, each $\lambda \in \cP$ can be written as in \eqref{eq31}. The restriction of $\lambda_{n+1}$ to $\lambda_n$
 is equivalent to
\begin{equation} \label{eq33}
\begin{array}{ccccccc}
a_n^1 & = & a_{n+1}^1 &+& 2 a_{n+1}^2 &+& a_{n+1}^3,\\
a_n^2 & = & a_{n+1}^3 &+& 2 a_{n+1}^4 &+& a_{n+1}^5,\\
\vdots &  &  \vdots &&  \vdots && \vdots\\
a_n^{2^n-1} & = & a_{n+1}^{2^{n+1}-3} & + & 2 a_{n+1}^{2^{n+1}-2} & + & a_{n+1}^{2^{n+1}-1}.
\end{array}
\end{equation}
Set $b_n := a_n^1 + \ldots + a_n^{2^n-1}$. Adding the equations in \eqref{eq33}
we obtain
\begin{equation} \label{eq34}
b_n = b_{n+1} + (a_{n+1}^2 + \ldots + a_{n+1}^{2^{n+1}-2}).
\end{equation}
Thus $b_1 \geq b_2 \geq b_3 \geq \ldots$ and, consequently, $b_{n_0} = b_{n_0+1} = b_{n_0 +2} + \dots$ for
some $n_0$.  Again \eqref{eq34} shows that $a_n^2 = \ldots = a_n^{2^n-2} = 0$ and $a_n^1 = a_{n+1}^1$, $a_n^{2^n-1} =
a_{n+1}^{2^{n+1}-1}$ for every $n > n_0$. Putting these facts together, we see that if $\lambda \in \cP^+$, then
$\lambda = \proj (a' \omega_n^1 + a'' \omega_n^{2^n-1})$. In particular, $\cP^+$ consists of the lattice points in a cone of
dimension two.

\noindent
(iii) In this case we are going to show that $\cP^+ = 0$.
Again, consider $\lambda = \proj \lambda_n \in \cP$. Let 
\begin{equation*} \label{eq35}
\lambda_n  =  \lambda_n^1 \vep_n^1 + \ldots + \lambda_n^{2^n-1} \vep_n^{2^n-1} 
=  a_n^1 \omega_n^1 + \ldots + a_n^{2^n-1} \omega_n^{2^n-1}.
\end{equation*}
Since $a_n^i = \lambda_n^i - \lambda_n^{i+1}$ for $1 \leq i \leq 2^n-2$ and $a_n^{2^n-1} = \lambda_n^{2^n-1}$, 
the fact that $\lambda_{n+1}$ restricts to $\lambda_n$ is equivalent to the equations
\begin{equation} \label{eq36}
\begin{array}{ccccccc}
a_n^1 & = & a_{n+1}^2 &+& 2 a_{n+1}^3 &+& a_{n+1}^4,\\
a_n^2 & = & a_{n+1}^4 &+& 2 a_{n+1}^5 &+& a_{n+1}^6,\\
\vdots &  &  \vdots &&  \vdots && \vdots\\
a_n^{2^n-2} & = & a_{n+1}^{2^{n+1}-4} & + & 2 a_{n+1}^{2^{n+1}-3} & + & a_{n+1}^{2^{n+1}-2},\\
a_n^{2^n-1} & = & a_{n+1}^{2^{n+1}-2} & + & 2 a_{n+1}^{2^{n+1}-1} .&  &
\end{array}
\end{equation}
Assume that $\cP^+ \neq 0$ and let $\lambda \in \cP^+$ be a nonzero weight. Choose $n_0$ so that
$\lambda_{n_0-1} \neq 0$
and set $b_k := a_{n_0+k}^{2^k} + a_{n_0+k}^{2^k +1} + \ldots + a_{n_0+k}^{2^{n+k}-1}$ for
$k \geq 0$. Adding the appropriate equations from \eqref{eq36} we obtain
\begin{equation} \label{eq37}
b_k = b_{k+1} + (a_{n_0 + k+1}^{2^{k+1} +1} + \ldots + a_{n_0+k+1}^{2^{n_0+k+1}-1}),
\end{equation}
which implies that $b_0 \geq b_1 \geq \ldots$. Hence there exists $k_0$ such that $b_{k_0} = b_{k_0+1} = \ldots$. We may assume 
that $k_0$ is the smallest such integer. 

If $k_0 = 0$, then \eqref{eq37} shows that
$$
a_{n_0 + 1}^{ 3} = \ldots = a_{n_0+1}^{2^{n_0+1}-1} = 0,
$$
which substituted in \eqref{eq36} implies 
$$
a_{n_0}^{ 2} = \ldots = a_{n_0}^{2^{n_0}-1} = 0,
$$
After another look at \eqref{eq36}, we obtain $\lambda_{n_0-1} = 0$, which contradicts the assumption that 
$\lambda_{n_0-1} \neq 0$.

If $k_0>1$,  then \eqref{eq37} shows that 
$$
a_{n_0 + k_0+1}^{2^{k_0+1} +1} = \ldots = a_{n_0+k_0+1}^{2^{n_0+k_0+1}-1} = 0,
$$
which substituted back into \eqref{eq36} gives
$$
a_{n_0 + k_0}^{2^{k_0} +1} = \ldots = a_{n_0+k_0}^{2^{n_0+k_0}-1} = 0.
$$
The last equation together with \eqref{eq37} implies $b_{k_0-1} = b_{k_0}$, which contradicts the choice of $k_0$. 
This proves that $\cP^+ = 0$.
\hfill $\square$
\end{example}

\vskip.5cm

Despite Example \ref{ex3} (iii), we can prove that $\cP^+ \neq 0$ under some natural assumptions. On the other hand, there are no
strictly dominant weights unless $G$ is root reductive.
 
 \vskip.5cm
 
\begin{proposition} \label{prop_P^+} $\phantom{x}$

\noindent
{\rm (i)} If $G$ is pure, then $\cP_B^+ \neq 0$ for any Borel subgroup $B \subset G$.

\noindent
{\rm (ii)} For any $G$ there exists a Borel subgroup $B$ such that $\cP_B^+ \neq 0$. 

\noindent
{\rm (iii)} Assume that $\cP_B^+$ contains a strictly dominant weight $\lambda$, i.e. such that $(\lambda_n, \alpha) > 0$
for every root $\alpha$ of $\gb_n$. Then $G$ is root reductive.  
\end{proposition}

\begin{proof}
(i) If $G$ is pure of type $B, C$, or $D$, then $\omega_{n+1}^1$ restricts to $\omega_n^1$ and thus
$\cP_B^+$ contains $\proj (a \omega_n^1)$ for every $a \in \Z_{\geq 0}$. If $G$ is pure of type $A$, then 
$\omega_{n+1}^1$ restricts to $\omega_n^1$ or $\omega_n^{r_n}$ and 
$\omega_{n+1}^{r_n}$ restricts to $\omega_n^1$ or $\omega_n^{r_n}$, which implies that every
$B_n$--dominant weight of the form $a' \omega_n^1 + a'' \omega_n^{r_n}$ extends to a $B_{n+1}$--dominant
weight. This shows that $\cP_B^+ \neq 0$. 

\noindent
(ii) Using induction we can construct compatible orders on the weights of the natural representation of $G_n$ in such a 
way that the maximal element among the weights of the natural representation of $G_{n+1}$ restricts to a weight
of the natural representation of $G_n$ and not to zero. Then, for every $a \in \Z_{> 0}$, $\proj (a \omega_n^1)$
is a nonzero element of $\cP_B^+$.

\noindent
(iii) Let $\alpha$ be a long root of $\gb_m$, and let $\alpha^1, \ldots, \alpha^{s_{m,n}}$ be the roots of $\gb_n$, $n >m$,
 which restrict to $\alpha$. Assuming that $\lambda_n$ is strictly dominant $\gg_n$--weight we conclude 
 that $(\lambda_n, \alpha^i) \geq 1/2$
for $1 \leq i \leq s_{m,n}$. This gives
$$
(\lambda_m, \alpha) = (\lambda_n, \alpha^1) + \ldots + (\lambda_n, \alpha^{s_{m,n}}) \geq \frac{1}{2}  s_{m,n},
$$
which is only possible if there is $n_0 >m$ so that $s_{m,n} = s_{m, n_0}$ for $n\geq n_0$. The latter condition implies that $G$ is root
reductive.
\end{proof}

\vskip.5cm

Every $\lambda \in \cP_B^+$ defines an irreducible $G$--module $V_{B}(\lambda)$ in the following way.
The weight $\lambda$ determines the direct system of highest weight modules $V_{B_n}(\lambda_n) \stackrel{e_n}{\to}
V_{B_{n+1}}(\lambda_{n+1})$, where $e_n$ maps the $B_n$--highest weight space of $V_{B_n}(\lambda_n)$ into the
$B_{n+1}$--highest weight space of $V_{B_{n+1}}(\lambda_{n+1})$. Then $V_B(\lambda)$ is defined as $\inj V_{B_n}(\lambda_n)$.

\vskip.5cm

\begin{example} \label{ex4}
 Let $G = SL(2^\infty)$ and let $B$ be any Borel subalgebra of $G$. Set 
$\lambda_n := \omega_n^1 + \omega_n^{2^n-1}$. The sequence $\{\lambda_n\}$ is a $B$--dominant weight of $G$
and hence the $G$--module $V_B(\lambda)$ is well--defined. Furthermore, in this case 
it is easy to check that $V_B(\lambda)$ is a weight module,
i.e.
$$
V_B(\lambda) = \oplus_\mu V_B(\lambda)^\mu, \quad {\text { where }} \quad
V_B(\lambda)^\mu = \{ v \in V_B(\lambda) \, | \, h \cdot v = \mu(h) v {\text { for every }} h \in \gh \}.
$$
This observation implies that, despite the fact that each of the modules $V_{B_n}(\lambda_n)$ 
is isomorphic to the adjoint representation of $G_n$, $V_B(\lambda)$ is not isomorphic to 
the adjoint representation of $G$ since the latter is not a weight module.
\hfill $\square$
\end{example}

\vskip.5cm

\section{The Weyl group $W_B$} \label{sec22} 
In this section we use the filtration
\eqref{eq12} to construct a group $W_B$ which plays the role that the Weyl group plays in
the classical Bott--Borel--Weil theorem.

First we consider the case when $G$ is not of type $B$. Let $W_n$ denote the Weyl group of $\gg_n$ and let 
$pr_n^i : \widetilde\gg_{n} = \gg_n^{\oplus s_{n}} \to \gg_n$ be the projection onto the $i^{th}$ direct summand for $1 \leq i \leq s_{n}$.
Then, for each pair $n, i$, the composition
\begin{equation} \label{eq41}
    \frak{g}_{n} \rightarrow \tilde\gg_{n} \overset{pr_n^{i}}{\rightarrow} {\frak{g}}_{n}^i\rightarrow \frak{g}_{n+1}
\end{equation}
is a root injection and hence yields an injective homomorphism of Weyl groups
$$
   \tau_n^i: W_{n} \to W_{n+1}.
$$
For every sequence $\{t_n\}_{n=1}^\infty$ with $1 \leq t_n \leq s_n$,
the injections $\tau_n^{t_n}$ form a direct system. Note that, if $\{t'_n\}$ and $\{t''_n\}$ are
two sequences which differ in finitely many positions only, then $\inj_{\tau_n^{t'_n}} W_n = \inj_{\tau_n^{t''_n}} W_n$.
We define an equivalence relation between sequences  by setting
$\{t'_n\} \sim \{t''_n\}$ if $t'_n = t''_n$ for large enough $n$, and denote the set of equivalence classes by $\cT$:
$$
\cT = \{\{t_n\} \, | \, 1 \leq t_n \leq s_n\}/\sim.
$$
The set $\cT$ consists of a single element if
$G$ is root reductive, and is uncountable otherwise. 
For any element $t \in \cT$ we  put $W^t := \inj_{\tau_n^{t_n}} W_n$, where $\{t_n\}$ is a representative of $t$. It is easy to
see that $W^t$ depends only on the type of $G$. Namely,
$W^t$ is isomorphic to $S_\infty$, the group of finite permutations of $\N$, if $G$ is of type A; to the group of 
signed finite permutations of $\N$ if $G$ is of type $C$; and to the group of 
signed finite permutations of $\N$ with even number of minus signs if $G$ is of type $D$.
Finally we put $W:=\dot{\times}_{t \in \cT}W^{t}$, where $\dot{\times}$ stands for restricted direct product. 

If $G$ is root reductive of type $B$ the definitions above still make sense. Moreover, $\cT$ consists of 
a single element $t$ and $W = W^t$ is isomorphic to the group of signed finite permutations of $\N$.

If $G$ is of type $B$ but is not root reductive, $\kappa_n$ is not a root injection for infinitely many $n$ 
and we need to modify the definitions above. 
Let $\mathring{W_n}$ denote the subgroup of $W_n$ generated by reflections along the long simple roots of
$\gb_n$.  It is clear that \eqref{eq41} maps root spaces corresponding to long roots of $\gg_n$ into
root spaces corresponding to long roots of $\gg_{n+1}$. Hence, again we have an injective homomorphism of groups
$$
   \tau_n^i: \mathring{W}_{n} \to \mathring{W}_{n+1}.
$$
We can now proceed as above to define $W^t :=\inj_{\tau_n^{t_n}} \mathring{W}_n$ and $W:=\dot{\times}_{t \in \cT}W^{t}$.
Note that $W^t$ is isomorphic to the group of finite permutations of $\N$. 

Next we define a {\it length function} $\ell_B: W \to \N \cup \{\infty\}$. 
Let $\ell_n$ denote the length function on $W_n$ determined by $\gb_n$.

\vskip.5cm

\begin{lemma} \label{le221}
For every $n \in \N$, every $1 \leq i \leq s_n$, and every $w \in W_{n}$ we have
$$ 
\ell_{n}(w) \leq \ell_{n + 1}( \tau_n^{i} (w)) .
$$
Furthermore, if $\ell_{n}(w) = \ell_{n + 1}( \tau_n^{i} (w))$ then for every reduced factorization $w = \sigma_1 \ldots \sigma_j$ into 
a product of simple reflections, $\tau_n^i(w) = \tau_n^i(\sigma_1) \ldots \tau_n^i(\sigma_j)$ is a  reduced factorization of 
$\tau_n^i(w)$. 
\end{lemma}

\begin{proof}
Let $\Delta_n = \Delta_n^+ \sqcup \Delta_n^-$ be the partition of the roots of $\gg_n$ into 
positive and negative corresponding to $\gb_n$, and let $\gamma: \Delta_n \to \Delta_{n+1}$
be the map corresponding to the injection \eqref{eq41}. Then $\gamma(\Delta_n^\pm) \subset \Delta_{n+1}^\pm$.
Set $\Phi_w := (w^{-1} \Delta_n^-) \cap \Delta_n^+$ and $\Phi_{\tau_n^i(w)}:= (\tau_n^i(w)^{-1} \Delta_{n+1}^-) \cap \Delta_{n+1}^+$.
The inclusion $\gamma(\Phi_w)  \subset \Phi_{\tau_n^i(w)}$ implies
$$
 \ell_{n}(w) = |\Phi_w| \leq |\Phi_{\tau_n^i(w)}| = \ell_{n + 1}( \tau_n^{i} (w)).
 $$
 Moreover, the equality $\ell_{n}(w) = \ell_{n + 1}( \tau_n^{i} (w))$ implies that $\gamma(\Phi_w)  = \Phi_{\tau_n^i(w)}$,
and thus $\gamma$ sends every simple root in $\Phi_w$  into a simple root of $\gb_{n+1}$. Consider a reduced factorization 
$w= \sigma_1 \ldots \sigma_j$ and let $\sigma_j$ be the reflection along the simple root $\alpha$ of $\gb_n$. Then 
$\tau_n^i(\sigma_j)$ is the reflection along the simple root $\gamma(\alpha)$ of $\gb_{n+1}$. Set $w' := \sigma_1 \ldots \sigma_{j-1}$.
Then we have $\ell_n(w') = \ell_{n+1}(\tau_n^i(w'))$, and we complete the proof by induction.
 \end{proof}

\vskip.5cm

\refle{le221} (and the observation that $\tau_n^i(\mathring{W}_n) \subset \mathring{W}_{n+1}$ if $G$ is of type $B$)
 implies that every element $w^t \in W^t$  has a 
well--defined, possibly infinite, {\it length} $\ell_B (w^t)$. We extend the definition of length to elements of $W$ by 
setting $\ell_B(w) := \sum_{t \in \cT} \ell_B(w^t)$ for $w = (w^t)_{t \in \cT}$. 
We now define $W_B$ as the subgroup 
of $W$ consisting of all elements $w \in W$ of finite length $\ell_B(w)$.
For $w = (w^t) \in W_B$ we say that the {\it support} of $w$ is $\{t^1, \ldots, t^l\} \subset \cT$ if 
$w^t = 1_{W^t}$ precisely when $t \not \in \{t^1, \ldots, t^l\}$.
Assume that $n_0$ is such that $w^{t^i} \in W_{n_0}$ for $1 \leq i \leq l$. Then $w^{t^i} \in W_{n}$ for $n \geq n_0$ and $1 \leq i \leq l$.
It is not necessarily true that $w^{t^1}, \ldots, w^{t^l}$
commute in $W_n$. If, however, the sequences $(t_{n_0}^i, \ldots, t_{n-1}^i)$ for $i = 1, \ldots, l$ are distinct, then
$w^{t^1}, \ldots, w^{t^l}$ commute in $W_n$, and define an element $w(n) := w^{t^1} \ldots w^{t^l} \in W_n$. Since
$t^1, \ldots, t^l$ are distinct, there exists $n_1$ such that the sequences $(t_{n_0}^i, \ldots, t_{n_1-1}^i)$ for $i = 1, \ldots, l$
are distinct and hence $w^{t^1}, \ldots, w^{t^l}$ commute in $W_n$ for every $n \geq n_1$. For the rest of the paper, whenever
for an element $w \in W_B$ we consider the elements $w(n)$, we will assume that $n \geq n_1$. 
If $G$ is of type $B$, then $w(n) \in \mathring{W}_n$. 

\vskip.5cm

\begin{example} \label{ex10}
For $G = SL(2^{ \infty })$ we have $W \cong \dot{\times}_\cT S_{\infty}$.
If $B$ is the upper triangular Borel subgroup of $G$, then $W_B = W$.
If, on the other hand, $B$ is the interlacing Borel subgroup, $W_B$ is trivial.
\hfill $\square$
\end{example}

\vskip.5cm

\begin{proposition}\label{prop221} 
If $W_B$ contains an element of length $l$, then $W_B$ contains elements of all lengths
from $0$ through $l$. $W_B$ may be finite or infinite and 
may or may not contain an element of maximal length. In addition, 
for  fixed $G$ and variable $B$, any non--negative integer can appear as a maximal 
possible length of an element in $W_B$.
\end{proposition}

\begin{proof}
The first statement follows from the generalization of Lemma \ref{le221} discussed above:
Let $w \in W$ with $\ell_B(w) = l$. For any  reduced expression for $w(n)$, 
any subword of $w(n)$ is well--defined and represents an element of $W_B$. Subwords of $w(n)$
will provide elements of $W$ of any length between $0$ and $l$. 

The remaining  statements are rather straightforward and we omit their proofs here.
\end{proof}

\vskip.5cm

It is not difficult to see that in general the group $W_B$ does not act on $\cP$, i.e.
there exist $w \in W_B$ and  $\lambda = \proj \lambda_n$ for which $\{w(n)(\lambda_n)\}$
is not an inverse system of weights. Here is a simple example.

\vskip.5cm

\begin{example} \label{ex11}
Let $B$ be the upper triangular Borel subgroup of $SL(2^\infty)$. Consider $w \in W_B$ given by
$$
w^t = \begin{cases} (12) & {\text { if }} t=(1, 1, \ldots)\\
1_{S_\infty} & {\text { otherwise}},
\end{cases}
$$
where the transposition $(12)$ is understood  as an element of the symmetric group $S_\infty$. 
Let $\lambda = \proj \lambda_n$ be a weight such that $a_n^1 = n$ and  $a_n^{2^{n-1}+1} = -1$ in the notation of 
Example \ref{ex3}(i). Then $w(n+1)(\lambda_{n+1}) = \lambda_{n+1} - (n+1) (\vep_n^1 - \vep_n^2)$
restricts to $\lambda_n - (n+1) (\vep_n^1 - \vep_n^2)$, while $w(n)(\lambda_{n}) = \lambda_{n} - (n) (\vep_n^1 - \vep_n^2)$,
which shows that $w(n)(\lambda_n)$ is not an inverse system of weights. 
\hfill $\square$
\end{example}

\vskip.5cm

Despite this example, we are going to show that if $\lambda \in \cP^+$ then $w(\lambda)$ is a well--defined 
element of $\cP$ for any $w \in W_B$. 
We will also define an analog of the "dot" action in the finite--dimensional case. 
 Recall that, for a finite--dimensional reductive Lie algebra $\gg'$ with 
fixed Cartan subalgebra $\hh'\subset\gg'$ and Borel subalgebra $\bb'\subset\gg'$, $\bb'\supset \hh'$, 
the dot action of a Weyl group element $w'$ on a weight $\mu'\in(\hh')^*$ is defined as 
$w'(\mu'+\rho_{\bb'})-\rho_{\bb'}$, where $\rho_{\bb'}$ is the half--sum of roots of $\bb'$. One writes 
$w'\cdot \mu' \eqdef w'(\mu'+\rho_{\bb'})-\rho_{\bb'}$.
In the case of the diagonal ind--group $G$, for any $\lambda = \proj \lambda_n$ and $w \in W_B$ it is
natural to consider the weights $\{w(n)(\lambda_n + \rho_n) - \rho_n\}$, where $\rho_n$ denotes
the half--sum of the roots of $\gb_n$. 

To prove  results about the action of $W_B$ on weights we need additional notation.
If $\alpha'$ is a root of $\gb_m$ and $\alpha''$ is a root of $\gb_n$ with $n \geq m$, we say that $\alpha''$ is
{\it a successor of $\alpha'$} if $\alpha'' \in \gh_n^*$ restricts to $\alpha' \in \gh_m^*$. If, in addition, $n = m+1$, we say
that $\alpha''$ is {\it an immediate successor} of $\alpha'$. Every root of $\gb_m$ has exactly $s_m$ immediate successors.
The set of successors $\cS^\alpha$ of a root $\alpha$ of $\gb_m$ has a natural structure of a directed tree --- 
every element is connected with
its immediate successors. If $\alpha$ is a root of $\gb_m$, then $\cS^\alpha = \sqcup_{n\geq m} \cS_n^\alpha$, where
$\cS_n^\alpha$ is the
set of successors of $\alpha$ of level $n$, i.e. those successors of $\alpha$ which are roots of $\gb_n$.
 Furthermore, given $\lambda \in \cP$, we assign integer labels to all nodes of this tree in a natural way:
the node $\alpha' \in \cS_n^\alpha$ is labeled by $\frac{2(\lambda_n, \alpha')}{(\alpha', \alpha')}$. 
It is clear that the sum of the labels
of the elements of $\cS_n^\alpha$ is the same for all $n$ and equals $\frac{2(\lambda_m, \alpha)}{(\alpha, \alpha)}$.

\vskip.5cm

\begin{proposition} \label{prop223} $\phantom{x}$

\noindent
{\rm (i)} If $w \in W_B$, then $w \cdot 0$ is  a well--defined element of $\cP$, i.e. $\{w(n) \cdot 0 = w(n)(\rho_n) - \rho_n\}$ is an inverse system 
of weights of $G$.

\noindent
{\rm (ii)} If $w \in W_B$ and $\lambda \in \cP_B^+$, then $w(\lambda)$ is  a well--defined element of $\cP$.

\noindent
{\rm (iii)} If $w \in W_B$ and $\lambda \in \cP_B^+$, then $w \cdot \lambda$ is a well--defined element of $\cP$.
\end{proposition}

\begin{proof}
(i) Since $w \in W_B$ we have $w^{-1} \in W_B$ as well and $w(n)^{-1} = w^{-1}(n)$. The proof of Lemma \ref{le221} applied
to $w^{-1}$ implies that the set $\Phi_{w(n+1)^{-1}}$ projects onto the set $\Phi_{w(n)^{-1}}$ and the 
 formulas, cf. \cite{DR},
$$
w(n) \cdot 0 = w(n)(\rho_n) - \rho_n = - \sum_{\alpha \in \Phi_{w(n)^{-1}}} \alpha
$$
and 
$$
w(n+1) \cdot 0 = w(n+1)(\rho_{n+1}) - \rho_{n+1} = - \sum_{\alpha \in \Phi_{w(n+1)^{-1}}} \alpha
$$
imply that $w(n+1) \cdot 0$ restricts to $w(n) \cdot 0$.

\noindent
(ii) Let $w = (w^t) \in W_B$ have support $t^1, \ldots, t^l$ and let 
$m$ be such that $w^{t^i}$ for $1 \leq i \leq l$ all belong to $W_n$ and commute
in $W_n$ for $n \geq m$.
Let $w(m) = \sigma_{\alpha^1} \ldots \sigma_{\alpha^q}$ be a reduced expression of $w(n)$. 
Then, for $n \geq m$, $w(n) = \sigma_{\alpha_n^1} \ldots \sigma_{\alpha_n^q}$ is a reduced expression of
$w(n)$ and $\alpha_n^i$ is a successor of $\alpha^i$ for $1 \leq i \leq q$. 
Furthermore, the sequence $\alpha^i = \alpha_m^i, \alpha_{m+1}^i, \alpha_{m+2}^i, \ldots$ is a path in $\cS^{\alpha^i}$. 
Since $\lambda$ is dominant, i.e. all labels in $\cS^{\alpha^i}$ corresponding to $\lambda$ are non--negative integers, 
there exists $n_0 \geq m$ such that the labels on each of the paths  
$\alpha^i = \alpha_m^i, \alpha_{m+1}^i, \alpha_{m+2}^i, \ldots$ of level $n \geq n_0$ stabilize. For $n \geq n_0$ we have
\begin{equation} \label{eq51}
w(n)(\lambda_n) = \lambda_n - \left( \sum_{1\leq i \leq q} \frac{2(\lambda_n, \alpha_n^i)}{(\alpha_n^i, \alpha_n^i)} \, \alpha_n^i - 
\sum_{1 \leq i < j \leq q} \frac{2(\lambda_n, \alpha_n^j)}{(\alpha_n^j, \alpha_n^j)} \frac{2(\alpha_n^j, \alpha_n^i)}{(\alpha_n^i, \alpha_n^i)}
\, \alpha_n^i + \ldots \right).
\end{equation}
Now consider the restriction of $w(n+1)(\lambda_{n+1})$ to $\gh_n^*$. By the definition of $\lambda$, $\lambda_{n+1}$ 
restricts to $\lambda_n$ and, by the stabilization of the labels along the paths 
$\alpha^i = \alpha_m^i, \alpha_{m+1}^i, \alpha_{m+2}^i, \ldots$,
$$ 
 \sum_{1\leq i \leq q} \frac{2(\lambda_{n+1}, \alpha_{n+1}^i)}{(\alpha_{n+1}^i, \alpha_{n+1}^i)} \, \alpha_{n+1}^i - 
\sum_{1 \leq i < j \leq q} \frac{2(\lambda_{n+1}, \alpha_{n+1}^j)}{(\alpha_{n+1}^j, \alpha_{n+1}^j)} 
\frac{2(\alpha_{n+1}^j, \alpha_{n+1}^i)}{(\alpha_{n+1}^i, \alpha_{n+1}^i)}
\, \alpha_{n+1}^i + \ldots 
$$
restricts to
$$
 \sum_{1\leq i \leq q} \frac{2(\lambda_n, \alpha_n^i)}{(\alpha_n^i, \alpha_n^i)} \, \alpha_n^i - 
\sum_{1 \leq i < j \leq q} \frac{2(\lambda_n, \alpha_n^j)}{(\alpha_n^j, \alpha_n^j)} 
\frac{2(\alpha_n^j, \alpha_n^i)}{(\alpha_n^i, \alpha_n^i)}
\, \alpha_n^i + \ldots.
$$
These observations together with \eqref{eq51} and its analog with $n+1$ in place of $n$ imply that 
$w(n+1)(\lambda_{n+1})$ restricts to $w(n)(\lambda_n)$. This completes the proof of (ii).

 \noindent
 (iii) The statement follows from (i), (ii), and the obvious formula
 $$
 w(n) \cdot \lambda_n = w(n)(\lambda_n) + w(n) \cdot 0. 
 $$
\end{proof}

\vskip.5cm

\begin{proposition} \label{prop41}
(i) Let $\sigma \in W_B$ be an element of length one and let $\lambda \in \cP$ be such that $\sigma(n) \cdot \lambda_n$ is
dominant for large enough $n$. Then there exists $n'$ such that, for $n \geq n'$, $\sigma(n) = \sigma_{\alpha_n}$, 
 $(\lambda_n, \alpha_n)$ does not depend on $n$, and $(\lambda_n, \alpha) = 0$ for
every successor $\alpha \in \cS_n^{\alpha_{n'}}$ different from $\alpha_n$;

\noindent
(ii) Let $w \in W_B$ and $\lambda \in \cP$ be such that $\mu_n := w(n) \cdot \lambda_n$ is $\gb_n$--dominant for large 
enough $n$. Then $\mu := \proj \mu_n$ is a well--defined element of $\cP_B^+$.
\end{proposition}

\begin{proof}
(i) Since $\ell_B(\sigma) = 1$, $\sigma(n) = \sigma_{\alpha_n}$ where $\alpha_{n+1}$ is an immediate successor of $\alpha_n$
for $n \geq n_0$ and both $\alpha_n$ and $\alpha_{n+1}$ are simple roots of the respective 
Borel subalgebras $\gb_{n}$ and $\gb_{n+1}$. 
The label of $\lambda$ at $\alpha_n$ in $\cS^{\alpha_{n_0}}$ is negative, while all other labels of $\lambda$ in $\cS^{\alpha_{n_0}}$ 
are non--negative. This implies that the labels of $\lambda$  along the path $\alpha_{n_0}, \alpha_{n_0 + 1}, \ldots$ are non--increasing.
Let $\beta_{n_0+1} = \alpha_{n_0+1} + \alpha_{n_0+1}'$ be a root of $\gb_{n_0 +1}$ higher than $\alpha_{n_0 +1}$
and let $\beta_n = \alpha_n + \alpha_n'$ be a successor of $\beta_{n_0+1}$. Note that $\beta_n$ is uniquely determined by 
$\alpha_n$. We have
$$
(\lambda_n,\beta_n) = (\lambda_n, \alpha_n) + (\lambda_n, \alpha_n'),
$$
which implies
$$
(\lambda_n, \alpha_n) \geq - (\lambda_n,\beta_n).
$$
Since $\{(\lambda_n,\beta_n)\}$ is a non--increasing sequence of non--negative integers or half--integers, we conclude that
the sequence $\{(\lambda_n, \alpha_n)\}$ is bounded, and hence it stabilizes. Noting that $(\lambda_{n+1}, \alpha_{n+1}) = (\lambda_n, \alpha_n)$
implies that $(\lambda_{n+1}, \alpha) = 0$ for every immediate successor of $\alpha_n$ other than $\alpha_{n+1}$ concludes the
proof of (i). 

\noindent
(ii) Write $w(n) = \sigma_{\alpha_n}^1 \ldots \sigma_{\alpha_n}^q$ as in the proof of Proposition \ref{prop223} (ii). 
As in (i) we prove that the labels along the paths $\{\alpha_n^i\}$ for $1 \leq i \leq q$ stabilize, and then repeat the 
argument in the proof of   Proposition \ref{prop223} (ii).
\end{proof}

\vskip.5cm

\section{$G/B$ and the Bott--Borel-Weil theorem} 

Recall that an {\it ind--variety} $X = \inj X_n$ is determined by a sequence of morphisms of algebraic varieties
$$
X_{1} \overset{\varphi_{1}}{\to} \ldots\to X_{n} \overset{\varphi_{n}}\to X_{n+1} \rightarrow   \ldots,
$$
see for instance \cite{Sh}, \cite{DPW}.
We denote by $\mathcal{O}_{X_{n}}$ the structure sheaf of $X_{n}$ and we define the 
{\it structure sheaf $\mathcal{O}$ of $X$} as the inverse limit 
$\lim \limits_{\leftarrow} \mathcal{O}_{X_{n}}$. More generally, a $sheaf$ $\mathcal{F}$ $on$ 
$X$ is by definition the limit of an inverse system of sheaves $\mathcal{F}_{n}$ on 
$X_{n}$, and $\mathcal{F}$ is a $sheaf$ $of$ $\mathcal{O}$--modules whenever 
$\{\mathcal{F}_{n}\}$ is an inverse system of sheaves of $\mathcal{O}_{X_{n}}$--modules. A sheaf 
of $\mathcal{O}$--modules is {\it locally free of rank $r$} whenever each $\mathcal{F}_{n}$ is 
locally free of rank $r$. In what follows we will also call a locally free sheaf of $\mathcal{O}$--modules 
a {\it vector bundle} on $X$.

Assume now that all $X_{n}$  are proper. Then it is well--known that the cohomology $H^{\cdot}(X, E)$ 
of any vector bundle $E = \lim \limits_{\leftarrow} E_{n}$ of finite--rank on $X$ is canonically isomorphic 
to the inverse limit $\lim \limits_{\leftarrow} H^{\cdot} (X_{n}, E_{n})$, see \cite{W}, \cite{DPW}.

In this paper we consider the ind--varieties $G/B$ and $G/P$, where $G=\lim \limits_{\rightarrow} G_{n}$ is 
a diagonal ind--group and $B = \lim \limits_{\rightarrow} B_{n}$ or $P=\lim \limits_{\rightarrow} P_{n}$ are 
respectively direct limits of Borel subgroups $B_{n}\subset G_{n}$ or parabolic subgroups $P_{n} \subset G_{n}$. 
More precisely, if $B=\lim \limits_{\rightarrow} B_{n}$ with $B_{n}=G_{n} \cap B_{n+1}$, or 
$P=\lim \limits_{\rightarrow} P_{n}$ with $P_{n}=G_{n} \cap P_{n+1}$, the embeddings 
$G_{n} \to G_{n+1}$ induce closed immersions $G_{n}/B_{n} \to G_{n+1}/B_{n+1}$ 
and $G_{n}/P_{n} \to G_{n+1}/P_{n+1}$ of proper smooth varieties. In what follows, we denote 
the corresponding ind--varieties by $G/B$ and $G/P$.

If $\lambda \in \cP$, the line bundles 
$(\mathcal{O}_{G_{n}/B_{n}})_{-\lambda_{n}}$ form an inverse system, and hence determine a line bundle (or a locally free 
sheaf of $\mathcal{O}$--modules of rank one) $\mathcal{O}_{-\lambda}$ on $G/B$. Recall that, by definition,
$( \mathcal{O}_{G_{n}/B_{n}} )_{-\lambda_{n}} $ is the $G_{n}$--equivariant line bundle on $G_{n}/B_{n}$ 
whose geometric fiber at the closed point $B_{n} \in G_{n}/B_{n}$ is the $B_{n}$--module $\Bbb{K}^{-\lambda_{n}}$.

More generally, if $E_{n}$ is a $G_{n}$--equivariant vector bundle (or, for short, {\it $G_n$--bundle}) on 
$G_{n}/P_{n}$, then the vector bundle $\displaystyle E=\lim_\gets E_n$ on $G_{n}/P_{n}$ is by definition 
$G$--equivariant, and each cohomology group $H^{j}(G/P, E)$ is a dual $G$--module, being an inverse 
limit of $G_{n}$--modules $H^{j}(G_{n}/P_{n}, E_{n})$.

The Bott--Borel--Weil theorem computes the cohomology 
$H^{\cdot}(G_{n}/B_{n}, (\mathcal{O}_{G_{n}/B_{n}})_{-\lambda_{n}})$ for each weight $\lambda_{n}$, see 
\cite{B}, \cite{D1}, \cite{D2}, and \cite{BW}. It is the following result.

\vskip.5cm

\begin{theorem}[Bott--Borel--Weil, \cite{BW}, \cite{B}] \label{th1}
If there exists a (necessarily unique) 
$w_{n} \in W_n$ such that $w_n \cdot \lambda_n$ is a $B_n$--dominant weight of $G_n$, then
$$
    H^{j}(G_{n}/B_{n},(\mathcal{O}_{G_{n}/B_{n}})_{-\lambda_{n}}) \cong \left\{
\begin{array}{ccc}
V_{B_{n}}(w_{n}\cdot \lambda_{n})^{*} & {\text{ for }} & j=\ell_{n}(w_{n}) \\
    \\
 0 & {\text{ for }} & j \neq \ell_{n}(w_{n}). \\
    \end{array} \right.
$$
If $w_n$ as above does not exists, then
$$
    H^{\cdot}(G_{n}/B_{n}, (\mathcal{O}_{G_{n}/B_{n}})_{-\lambda_{n}})=0.
$$
\end{theorem}

\vskip.5cm

An immediate corollary of \refth{th1} is that, for a fixed $\lambda \in \cP$, 
there is at most one $j$ for which the cohomology group $H^{j}(G/B, \mathcal{O}_{-\lambda})$ can be nonzero. 
This follows from \refth{th1} and from the fact that 
$H^{j}(G/B, \mathcal{O}_{-\lambda})=\lim \limits_{\leftarrow} H^{j}(G_{n}/B_{n}, (\mathcal{O}_{G_{n}/B_{n}})_{-\lambda_{n}}).$
The following theorem provides a much stronger statement. It is an analog of the Bott--Borel--Weil theorem and is the central result
in this paper.


\vskip.5cm

\begin{theorem} \label{th2} 
Let $G$ be a diagonal ind--group, let $B$ be a Borel subgroup of $G$, and  let $\lambda \in \cP$.
Then $H^j(G/B, \cO_{-\lambda}) \neq 0$ for at most one value of $j$. More precisely,
$H^j(G/B, \cO_{-\lambda}) \neq 0$ if and only if there exists $w \in W_B$ such that $w \cdot \lambda \in \cP_B^+$.
In the latter case we have an isomorphism of dual $G$--modules
$$
H^j(G/B, \cO_{-\lambda}) \cong  V_B(w \cdot \lambda)^*.
$$
  \end{theorem}
  
  \vskip.5cm

Before we prove \refth{th2} we state two results necessary for the proof. 
Let $G' \subset G''$ be reductive algebraic groups with Lie algebras $\gg' \subset \gg''$ respectively. 
Assume that $B'' \subset G''$ is a Borel subgroup of $G''$ and that $B' := G' \cap B''$ is a Borel subgroup
of $G'$. Then we have a close immersion of homogeneous spaces $G'/B' \to G''/B''$. Denote the Weyl groups of $G'$ and $G''$ by
$W'$ and $W''$ respectively. If $G'$ (respectively, $G''$) is of type $B$ or product of groups of type $B$, 
denote by $\mathring{W}'$ (respectively, 
$\mathring{W}''$) the subgroup of $W'$ (respectively, $W''$) generated by reflections along the long simple roots of the
corresponding Borel subalgebra.
Let $\lambda''$ be a weight of $B''$ which restricts to the weight $\lambda'$ of $G'$. 
Assume that there exist $w' \in W'$ and $w'' \in W''$ both of length $j$ and 
such that $w'' \cdot \lambda''$ and $w' \cdot \lambda'$ are dominant weights. The natural map
\begin{equation} \label{eq51}
H^j(G''/B'', \cO_{-\lambda''}) \to H^j(G'/B', \cO_{-\lambda'})
\end{equation}
is a homomorphism of nontrivial $G'$--modules. 

\vskip.5cm

\begin{proposition}[Tsanov, \cite{T}]   \label{thDotAction} In the notation above the following statements hold.

\noindent
{\rm (i)}  Assume that $\gg'$ is a root subalgebra of $\gg''$ and consider $W'$ as a subgroup of $W''$.
Then \eqref{eq51} is nonzero if and only if $w'' = w' \in W'$.

\noindent
{\rm (ii)} Assume that $\gg' \cong B_r \oplus B_r$ and $\gg'' \cong D_{2r+1}$ as in section \ref{sec21}. If  $w'' = w' \in \mathring{W}'$ then
\eqref{eq51}  is nonzero.
\end{proposition}

\vskip.5cm
\begin{corollary} \label{cor51}
Let $G' = G_1 \times \ldots \times G_s \cong (B_{r'})^s$ with $s>1$, let $G'' \cong B_{r''}$ and assume that 
the embedding $\kappa: G' \to G''$
is analogous to $\kappa_n$ from \eqref{eq13}. In other words, the natural representation $V''$ of $G''$ decomposes
as $V'_1 \oplus \ldots \oplus V'_s \oplus \K^z$, where each $V'_i$ is isomorphic to the natural representation of $B_{r'}$.
Consider the diagram
\begin{equation} \label{eq78}
\xymatrix{
\gg'  \ar[rr]^{\kappa} \ar[dr]_{\psi} & &\gg'' \\
& \bar{\gg}' \ar[ur]_{\theta} &},
\end{equation}
where $\bar{\gg}'$ is defined analogously to $\bar{\gg}_n$ from \eqref{eq13}. 

\noindent
{\rm (i)} If \eqref{eq51} is nonzero then $w'' \in \mathring{W}''$.

\noindent
{\rm(ii)} If $w' \in \mathring{W}'$, then \eqref{eq51} is nonzero if and only if $w'' = w'$.
\end{corollary}

\vskip.5cm
\begin{proof} The second statement follows from Proposition \ref{thDotAction}. Here is the proof of (i).
Denote by $\bar{G}'$ the simply--connected algebraic group with Lie algebra $\bar{\gg}'$. Assume that
 \eqref{eq51} is nonzero. The fact that \eqref{eq51} is nonzero
  implies that the map
$$
H^j(G''/B'', \cO_{-\lambda''}) \to H^j(\bar{G}'/\bar{B}', \cO_{-\bar{\lambda}'})
$$
is nonzero, where $\bar{\lambda}'$ is the restriction of $\lambda''$ to $\bar{G}'$. Since $\theta$ is a root injection,
Proposition \ref{thDotAction}(i) implies that $w'' \in \mathring{W}''$ if $s$ is even. If $s$ is odd, 
Proposition \ref{thDotAction}(i) implies that $w''$ is contained in the subgroup 
$^1\mathring{W}''$ of $W''$ generated by reflections along 
long simple roots corresponding to the components of $\bar{G}'$ of type $D$ and by reflections along the simple roots
corresponding to the component of $\bar{G}'$ of type $B$. We can use a different way of combining the $G'$--constituents
$V'_1, \ldots, V'_s$ 
of the natural representation of $G''$ to obtain a diagram analogous to \eqref{eq78} but with different $\bar{\gg}'$ and
different maps $\kappa$ and $\theta$. Repeating the argument above we conclude that $w'' \in {^2\mathring{W}''}$ for
a subgroup $^2\mathring{W}''$ analogous to $^1\mathring{W}''$. Note that $^1\mathring{W}'' \cap {^2\mathring{W}''} =
\mathring{W}''$ as long as we choose a different component of type $B$ of $\bar{\gg}'$. This completes the proof. 
\end{proof}

\vskip.5cm

\begin{proposition}[Corollary 5.4.1, \cite{DR}] \label{prop_DR}
Let $\gg''$ equal the direct sum of $s$ isomorphic copies $\gg''_1, \ldots, \gg''_s$ of $\gg'$
so that $\gg'$ projects isomorphically onto each subalgebra $\gg''_i$.
Assume that
$$
H^j(G''/B'', \cO_{-\lambda''}) = H^j(G_1''/B_1'', \cO_{-\lambda_1''}) \otimes H^0(G_2''/B_2'', \cO_{-\lambda_2''}) \otimes \ldots \otimes
H^0(G_s''/B_s'', \cO_{-\lambda_s''}),
$$
where $G'' = G_1'' \times \ldots \times G_s'' \cong (G')^s$ and $\lambda_1'', \ldots, \lambda_s''$ are the restrictions of $\lambda''$ to
$\gg_1'', \ldots, \gg_s''$. Then \eqref{eq51} is a nonzero homomorphism.
\end{proposition}

\vskip.5cm

\begin{proof}[Proof of Theorem \ref{th2}] For $n > m$ diagram \eqref{eq12} induces the commutative diagram
\begin{equation} \label{eq52}
\xymatrix{
\gg_m  \ar[rr]^{\delta_{m,n}} \ar[dr]_{\varphi_{m,n}}  & &\gg_{n} \\
& {\gg}_{m}^{n}
\ar[ur]_{\kappa_{m,n}} &},
\end{equation}
where $\delta_{m,n} = \delta_{n-1} \circ \ldots \circ \delta_m$
 and $\varphi_{m,n}$ and $\kappa_{m,n}$ are defined in the 
obvious way. By definition $\tilde{\gg}_m := \gg_{m}^{m+1}$ and ${\gg}_{m}^{n} \cong \gg_m^{\oplus s_m^n}$, where 
$s_m^n := s_m \ldots s_{n-1}$.
Given $\lambda \in \cP_B$, we denote the restriction of $\lambda$ to ${\gg}_{m}^{n}$ by $\lambda_m^n$.
Furthermore, for $n > m > k$ there exist maps $\delta_{k,m,n}$, $\varphi_{k}^{m,n}$, and $\kappa_{k,m}^{n}$ 
such that the diagram 
\begin{equation} \label{eq53}
\xymatrix{
\gg_k  \ar[rrr]^{\delta_{k,m}} \ar[drr]^{\varphi_{k,m}} \ar[dddrrr]_{\varphi_{k,n}} & & &\gg_{m} \ar[rrr]^{\delta_{m,n}} \ar[dr]^{\varphi_{m,n}} & & & \gg_n\\
& &{\gg}_{k}^{m} \ar[ur]^{\kappa_{k,m}} \ar[rr]^{\delta_{k,m,n}} \ar[rdd]^{\varphi_{k}^{m,n}} & & {\gg}_{m}^{n} \ar[urr]^{\kappa_{m,n}} & & \\
& & & & && \\
& && {\gg}_{k}^{n} \ar[uur]^{\kappa_{k,m}^{n}} \ar[uuurrr]_{\kappa_{k,n}} & &&}
\end{equation}
commutes. To simplify notation we set 
$$
\varphi_m^n := \left\{ \begin{array}{ccc}
\varphi_m & {\text { if }} n = m \\
&&\\
 \varphi_m^{n, n+1} & {\text { if }} n >m, \end{array} \right.
 $$
 
 $$
\kappa_m^n := \left\{ \begin{array}{ccc}
\kappa_m & {\text { if }} n = m+1 \\
&&\\
 \kappa_{m, m+1}^{n} & {\text { if }} n >m +1, \end{array} \right.
 $$
and $\gg_n^n := \gg_n$.
\vskip.5cm

Diagram \eqref{eq53} gives rise to a commutative diagram
\begin{equation} \label{eq54}
\begin{array}{ccccccccccccccc}
\ldots               & \gg_m&               & \to &                & \gg_{m+1} &                & \to &                  & \gg_{m+2}  &              & \to &                 & \gg_{m+3}          &\ldots\\
                               &     &\searrow &       &\nearrow &       & \searrow &       & \nearrow &       & \searrow &      & \nearrow &                   &\\
	        &  \ldots   &                & \gg_m^{m+1} &                & &                 & \gg_{m+1}^{m+2} &                 &  &                & \gg_{m+2}^{m+3}         &        &\ldots      & \\
	                        &     & &       &\searrow &       & \nearrow &       & \searrow &       & \nearrow &      &           &       &\\
	           &     &           & \ldots &                & \gg_m^{m+2} &                 &  &                 & \gg_{m+1}^{m+3} &                & \ldots &                &       & \\
	        	                        &     & &       & &       & \searrow &       & \nearrow &       &  &      &  &                   &\\
	            &     &           &  &                & \ldots &                 & \gg_m^{m+3} &                 & \ldots &                &  &                &        & \\
	            &&&&&&&&&&&&&&\\
	           &     &           &  &                &  &                 & \ldots &                 &  &                &  &                &  &      
\end{array}
\end{equation}
\vskip.5cm

The idea of the proof is to study the maps between the cohomology groups \\
 $H^j(G_m^n/B_m^n, \cO_{-\lambda_m^n})$ induced from \eqref{eq54}, 
where $G_m^n$ is the simply--connected  classical group with Lie algebra $\gg_m^n$ and $B_m^n = B \cap G_m^n$. More precisely, every injection 
in the direct system
\begin{equation} \label{eq57}
\gg_k \to \gg_k^{k+1} \to \gg_k^{k+2} \to \ldots
\end{equation}
splits into the direct product of embeddings to which Proposition \ref{prop_DR} applies. Similarly, 
Proposition \ref{thDotAction} and Corollary 
\ref{cor51} apply to every injection in the sequence
\begin{equation} \label{eq58}
\gg_m^n \to \gg_{m+1}^n \to \ldots \to \gg_{n-1}^n \to \gg_n.
\end{equation}

Assume first that $H^j(G/B, \cO_{-\lambda}) \neq 0$. Then there exists $k$ such that all 
maps between the cohomology groups $H^j(G_m^n/B_m^n, \cO_{-\lambda_m^n})$ induced from \eqref{eq54} for $n \geq m \geq k$ are
nonzero. Note that
$$
H^j (G_k^m/B_k^m, \cO_{-\lambda_k^m}) = \otimes_{t \in \cT_k^m} H^{j_t}((G_k^m)_t/(B_k^m)_t, \cO_{-(\lambda_k^m)_t}),
$$
where $\cT_k^m = \{(t_k, \ldots, t_{m-1}) \, | \, 1 \leq t_i \leq s_i\}$, $(G_k^m)_t$ and $(B_k^m)_t$ are the respective 
constituents of $G_k^m$ and $B_k^m$ corresponding to $t$, and $(\lambda_k^m)_t$ is the restriction of $\lambda$ to $(G_k^m)_t$. 
Let $w^t$ denote the element of the Weyl group  $(W_k^m)_t$ of $(G_k^m)_t$ such that $w^t \cdot (\lambda_k^m)_t$ is dominant.
We say that $t'' \in \cT_k^{m+1}$ is an {\it immediate successor of} $t' \in \cT_k^m$  if
$t' = (t_k, \ldots, t_{m-1})$ and $t'' = (t_k, \ldots, t_{m-1}, t_m)$; we denote this relation by $t' \prec t''$.
Kunneth's formula implies that $j_t = \sum_{t \prec t'} j_{t'}$. Hence there is a finite collection of sequences 
$$
t^i = (t_k^i, t_{k+1}^i, \ldots), {\text { for }} 1 \leq i \leq l
$$
such that, for $t \in \cT_k^m$,  $j_t \neq 0$ if and only if $t = (t_k^i, \ldots, t_{m-1}^i)$ for some $1 \leq i \leq l$.
Fix $m>k$ such that, for $t \in \cT_k^m$ with $j_t \neq 0$, $j_{t'} \neq 0$ for exactly one immediate successor $t'$ of $t$.
In particular $j_{t^i}$ stabilizes for $n >m$.
Noting that $(\lambda_k^{m+1})_{t''}$ is dominant for every $t \prec t'' \neq t'$ we conclude that $\tau_n^{t^i_{n}}(w^{t^i}) = w^{t^i}$
for $n > m$. The last equation means that we have well--defined elements $w^{t^i} \in W^{t^i}$ if $G$ is not of type $B$. 
If $G$ is of type $B$, Corollary \ref{cor51}(i) ensures that $w^{t^i} \in W^{t^i}$ if we repeat the argument above with $k+1$ in place 
of $k$.
Furthermore,
$\ell_B(w^{t^i}) = j_{t^i}$ and $w^{t^1}, \ldots, w^{t^l}$ define an element $w$ of $W_B$ of length $j$. 
The fact that 
$w \cdot \lambda \in \cP_B^+$ follows from Proposition \ref{prop41}. The existence of an isomorphism
$H^j(G/B, \cO_{-\lambda}) \cong  V_B(w \cdot \lambda)^*$ is obvious.

Conversely, assume that $w \in W_B$ satisfies $w \cdot \lambda \in \cP_B^+$. We need to show that there exists $k$ such that
all maps  
between cohomology groups $H^j(G_m^n/B_m^n, \cO_{-\lambda_m^n})$ corresponding to \eqref{eq54} 
with $n \geq m > k$ are nonzero. Assume that the support of $w$ is $t^1, \ldots, t^l$ and choose $k$
so that the sequences $t_k^i, t_{k+1}^i, \ldots$ for $1 \leq i \leq l$ are distinct. The fact  that 
$H^j(G_m^n, \cO_{-\lambda_m^n}) \to H^j(G_m^{n-1}, \cO_{-\lambda_m^{n-1}})$ is nonzero follows from
Propodition \ref{prop_DR}, while the fact that 
$H^j(G_m^n, \cO_{-\lambda_m^n}) \to H^j(G_{m-1}^{n}, \cO_{-\lambda_{m-1}^{n}})$ is nonzero follows 
from Proposition \ref{thDotAction} and Corollary \ref{cor51}. 
Finally, $H^j(G_{n+1}/B_{n+1}, \cO_{-\lambda_{n+1}}) \to H^j(G_{n}/B_{n}, \cO_{-\lambda_{n}})$
being the composition of 
$$
H^j(G_{n+1}^{n+1}, \cO_{-\lambda_{n+1}^{n+1}}) \to H^j(G_{n}^{n+1}, \cO_{-\lambda_{n}^{n+1}}) \quad {\text { and }} \quad
H^j(G_{n}^{n+1}, \cO_{-\lambda_{n}^{n+1}}) \to H^j(G_{n}^{n}, \cO_{-\lambda_{n}^{n}})
$$ 
is  nonzero.
\end{proof}

\vskip.5cm

\begin{example} \label{ex45}
 Let $G=SL(2^\infty)$. If $B$ is the upper triangular Borel subgroup, then \refth{th2}, together with 
 the explicit description of $W_{B}$ given above, implies that for each $j$ there are line 
bundles $\OO_{\lambda}$ with $H^j(G/B,\OO_{-\lambda})\neq 0$. If $B$ is the interlacing Borel subgroup of 
$SL(2^\infty)$, then $W_B$ is the trivial subgroup of $W$, and hence, by \refth{th2}, $H^j(G/B,\OO_{-\lambda})=0$ 
for all $B$--weights $\lambda$ and all $j> 0$. Moreover, in this case $H^j(G/B,\LL)=0$ for all $j>0$ and for any line 
bundle $\LL$ on $G/B$, as it is easy to show that any $\LL$ is $G$--equivariant, i.e. $\LL\simeq\OO_{-\lambda}$ 
for some $B$--weight $\lambda$. This implies that the above two homogeneous ind--spaces are not isomorphic as 
ind--varieties, and in particular that the interlacing Borel subgroup is not conjugate to the upper triangular Borel 
subgroup by an automorphism of $SL(2^\infty)$. \hfill $\square$
\end{example}
 
\vskip.5cm
 
Note that the group $W_B$ which we use in Theorem \ref{th2} is different from the Weyl group $W_F$ defined in \cite{NRW}
unless $G$ is root reductive. (In fact $W_F$ is a trivial group if $G$ is diagonal but not root reductive.) Nevertheless, if
$H^j(G/B, \cO_{-\lambda}) \neq 0$ for some $\lambda$, the $B$--weight $\lambda$ is cohomologically finite in the sense of 
\cite{NRW}. A question we do not answer in the present paper is whether the cohomological finiteness of $\lambda$ is sufficient
for $H^j(G/B, \cO_{-\lambda})$ to be nonzero.

\vskip.5cm

\section{$G/P$ and projectivity} \label{sec5}
Let $P$ be a parabolic subgroup of $G$ and let $B \subset P$ be a Borel subgroup of $G$.
It is easy to see that every  finite--dimensional simple $P$--module admits a $B$--highest weight, 
i.e. is the limit of the direct system of simple highest weight $P_n$--modules  for an inverse system of
$G_n$--weights $\lambda_n$. 
Setting $\lambda = \proj \lambda_n$, we denote by 
$M_\lambda$ the simple $P$--module with highest weight $\lambda$.  

\vskip.5cm

\begin{proposition} \label{prop51} Set $\cO(M_\lambda^*) := \proj (\cO_{G_n/P_n})(M_\lambda^*)$, where 
$(\cO_{G_n/P_n})(M_\lambda^*)$ is the usual $G_n$--equivariant bundle on $G_n/P_n$ with fibre $M_\lambda^*$.  
Then 
$H^j(G/P, \cO(M_\lambda^*)) \neq 0$ if and only if \\
$H^j(G/B, \cO_{-\lambda}) \neq 0$, and in that
case
$$
H^j(G/P, \cO(M_\lambda^*)) = H^j(G/B, \cO_{-\lambda}) \cong V_B(w \cdot \lambda)^*,
$$
where $w \in W_B$ and $w \cdot \lambda \in \cP_B^+$.
\end{proposition}

\begin{proof}
It is easy to see that $\cO(M_\lambda^*) \cong pr_* \cO_{-\lambda}$, $pr: G/B \to G/P$ being the natural submersion.
Moreover, the fibre of $pr$ equals $P/B = \inj P_n/B_n$, hence the classical Bott--Borel--Weil theorem implies $R^i pr_* \cO_{G/P} =0$
for $i >0$ and $pr_* \cO_{G/P} = \cO_{G/B}$.  This is sufficient to conclude that
$$
H^j(G/P, \cO(M_\lambda^*)) = H^j(G/B, \cO_{-\lambda}) 
$$
for any $j \in \Z_{\geq 0}$. The isomorphism $H^j(G/B, \cO_{-\lambda}) \cong V_B(w \cdot \lambda)^*$ is established in Theorem
\ref{th2}. 
\end{proof}

\vskip.5cm

We conclude this paper by discussing the projectivity of the ind--varieties $G/B$ and $G/P$.  
Recall that an ind--variety $X$ is projective, i.e. admits an embedding in the projective ind--space $\P^\infty$, if and only if
it admits a very ample line bundle $\cL$. 
An explicit criterion for the projectivity of $G/B$ (and, more generally, of $G/P$) 
when $G$ is root reductive is proved in \cite{DPW}. 

For diagonal ind--groups we have the following. 

\vskip.5cm

\begin{corollary} 
Let $G$ a diagonal ind--group and $B$ be a Borel subgroup of $G$. Then if $G/B$ is projective, $G$ is necessarily root reductive.
\end{corollary}

\begin{proof} If $\iota: G/B \to \P^\infty$ is a closed immersion then $\cL:= \iota^*(\cO_{\P^\infty}(1))$ is a very ample line bundle on $G/B$.
In other words, $\cL_{|G_n/B_n}$ is very ample for each $n$. Since $G_n$ is simply--connected for each $n$, 
$\cL_{|G_n/B_n} \cong (\cO_{G_n/B_n})_{-\lambda_n}$ for some strictly dominant weight $\lambda_n$ of $G_n$. The weights $\lambda_n$
form an inverse system and hence define a strictly dominant weight $\lambda = \proj \lambda_n$ of $G$. By Proposition \ref{prop_P^+} (iii)
$G$ is root reductive.
\end{proof}

\vskip.5cm

The following example shows that $G/P$ may be projective even if $G$ is not root reductive.

\vskip.5cm

\begin{example} \label{projective} 
Let $G = SL(2^\infty)$ and $P_n$ be the stabilizer of the span of the first $i$ standard basis vectors in $\C^{2^n}$. 
Then $\inj P_n$ is a well--defined maximal parabolic subgroup of $G$, and it is easy to see that  
 $G/P$ is isomorphic to the ind--Grassmannian  
of $i$--dimensional subspaces of $\C^\infty$. The latter is clearly projective. 
\end{example}


\end{document}